\documentclass[11pt,letterpaper]{article}
\usepackage[margin=1in]{geometry}
\usepackage[T1]{fontenc}
\usepackage{lmodern,mathtools,amssymb,amsthm,bbm,microtype,setspace}
\usepackage[round,authoryear]{natbib}
\usepackage{xcolor}
\definecolor{linkblue}{RGB}{22,67,114}
\definecolor{citegreen}{RGB}{35,102,78}
\definecolor{urlblue}{RGB}{0,91,150}
\usepackage[colorlinks,linkcolor=linkblue,citecolor=citegreen,urlcolor=urlblue]{hyperref}
\hypersetup{pdftitle={Weak and Strong Moments of Log-Concave Measures},
  pdfauthor={Haozhou Zhang},
  pdfsubject={A moment comparison for seminorms of log-concave random vectors},
  pdfkeywords={log-concave measure, weak moment, seminorm, relative entropy, mutual information}}
\newtheorem{theorem}{Theorem}[section]
\newtheorem{lemma}[theorem]{Lemma}
\newtheorem{proposition}[theorem]{Proposition}
\newtheorem{corollary}[theorem]{Corollary}
\theoremstyle{remark}
\newtheorem*{remark}{Remark}
\numberwithin{equation}{section}
\newcommand{\R}{\mathbb R}
\newcommand{\E}{\mathbb E}
\newcommand{\Prob}{\mathbb P}

\newcommand{\law}{\mathcal L}
\newcommand{\dd}{\,\mathrm d}
\newcommand{\KL}[2]{D(#1\mathbin\Vert #2)}
\newcommand{\bb}{\mathcal B_N}

\title{Weak and Strong Moments of Log-Concave Measures}
\author{Haozhou Zhang\\
Department of Mathematics and Statistics\\
University of Idaho\\
\href{mailto:zhan3063@vandals.uidaho.edu}{\texttt{zhan3063@vandals.uidaho.edu}}}
\date{}
\makeatletter
\renewcommand{\@maketitle}{%
  \begin{center}
    {\LARGE\@title\par}\vskip .7em
    {\begin{tabular}{c}\@author\end{tabular}\par}
  \end{center}}
\makeatother

\begin{document}
\maketitle
\begin{abstract}
Let \(X\in\R^d\) be a log-concave random vector, and let \(N\) be a seminorm
on \(\R^d\). For \(p\ge1\), we show that \(\|N(X)\|_p\) is at most
\(\E N(X)+C s_p\), where \(s_p\) is the supremum of
\(\|u(X)\|_p\) over all linear functionals \(u\) satisfying
\(|u(x)|\le N(x)\) for every \(x\in\R^d\). A choice for \(C\) is \(2+2e\).
This bound implies Lata{\l}a's conjecture on weak and strong moments.
\end{abstract}


\setlength{\abovedisplayskip}{8pt plus 2pt minus 1pt}
\setlength{\belowdisplayskip}{8pt plus 2pt minus 1pt}
\section{Introduction}\label{sec:intro}
Let \(X\in\R^d\) have law \(\mu\), and let \(N\) be a seminorm on \(\R^d\)
with dual unit ball
\[
 \bb=\{u\in(\R^d)^*: |u(x)|\le N(x)\text{ for all }x\in\R^d\}.
\]
Set
\[
 m=\E N(X),\qquad s_p=\sup_{u\in\bb}\|u(X)\|_p \quad (1\le p\le\infty).
\]
For log-concave \(X\) and a norm \(N\), \citet[Conjecture~1.1]{latala} conjectured that
\(\|N(X)\|_p\le C_1m+C_2s_p\) for \(1\le p<\infty\), with absolute constants \(C_1,C_2\).
In the Euclidean case, \citet{paouris} proved the bound and \citet{adamczak}
asked whether one can take \(C_1=1\).

\begin{theorem}\label{thm:main}
Let \(X\) be a log-concave random vector in \(\R^d\) and \(N\) a seminorm on \(\R^d\).
For \(1\le p\le\infty\),
\begin{equation}\label{eq:main}
 \|N(X)\|_p\le \E N(X-\E X)+(1+2e)s_p
 \le m+C_0s_p,\qquad C_0=2+2e.
\end{equation}
\end{theorem}
\noindent Here log-concavity means that for compact \(A,B\subset\R^d\) of positive
\(\mu\)-measure and \(0<\theta<1\),
\[
 \mu\big(\theta A+(1-\theta)B\big)\ge\mu(A)^\theta\mu(B)^{1-\theta}.
\]

Intuitively, Theorem~\ref{thm:main} states that whenever the \(L^p\) norm of
\(N(X)\) is much larger than \(\E N(X)\), some one-dimensional shadow \(u(X)\)
with \(u\in\bb\) must have a comparable \(L^p\) norm, regardless of the dimension.

\citet{ls} proved the comparison for \(\ell_r\) norms, \(1\le r<\infty\),
with constants of order at most \(r\). For arbitrary norms, the bound with
\(C_1=1\) was known for symmetric product log-concave laws, uniform laws on
\(\ell_r\) balls, and rotationally invariant log-concave laws
\citep[Theorem~6 and Proposition~13]{survey}.
\citet[Corollary~2.5]{sst} obtained the product case for \(p\ge2\) assuming
independent symmetric coordinates with log-concave tails.
\citet[Theorem~1.3]{hl} proved a related two-vector comparison for \(p\ge2\), with a
dimension-dependent coefficient on the mean in the isotropic case. For symmetric
laws, \citet[Conjecture~3]{lw} also asked for a bound on fluctuations about a lower moment.

The main difficulty of the proof lies in bounding the moments of maxima over finite subsets of \(\bb\)
uniformly in their cardinality. The expectation of each such maximum under \(\mu\)
is at most \(m\). To use such bounds for higher moments, we pass to
exponential moments and apply Gibbs duality, which reduces the problem to
bounding expected maxima under changed laws \(\nu\) minus a relative entropy penalty.
Fixing the distribution of a maximising index under \(\nu\), we split the
expected maximum into the best payoff at a positive mutual-information budget
and the remaining gain. An entropy estimate bounds the first term via
directional exponential moments, while entropy interpolation and label resampling
compare the remaining gain with the gain over independent selection under
\(\mu\), which is bounded by \(\E N(X-\E X)\) after centring. Choosing the budget
yields the Laplace comparison in Proposition~\ref{prop:laplace}, and a one-dimensional estimate then
bounds the directional transforms using \(s_p\), establishing Theorem~\ref{thm:main}.

\section{Entropy interpolation and the Laplace comparison}\label{sec:moments}

For probability laws \(P,Q\), let \(\KL P Q=\int\log(\mathrm dP/\mathrm dQ)\dd P\)
when \(P\ll Q\), and \(+\infty\) otherwise.

To compare selection rules, we will interpolate their conditional laws using
the following compactly supported case of displacement convexity
\citep[Proposition~2.3]{asz}.

\begin{lemma}\label{lem:interp}
Let \(\mu\) be compactly supported and log-concave, and let \(P_0,P_1\)
have finite entropy relative to \(\mu\). For \(0<\theta<1\) there is a
coupling \(X_j\sim P_j\) such that \(Z=(1-\theta)X_0+\theta X_1\) satisfies
\begin{equation}\label{eq:interp}
 \KL{\law(Z)}{\mu}\le (1-\theta)\KL{P_0}{\mu}+\theta\KL{P_1}{\mu}.
\end{equation}
\end{lemma}

\begin{proof}
Write \(\alpha_0=1-\theta\) and \(\alpha_1=\theta\). Suppose we match compact
sets \(A_j\) of \(P_j\)-mass \(w>0\). Interpolating points from these sets gives
points in \(B=\alpha_0A_0+\alpha_1A_1\), and log-concavity gives
\[
 w\log\frac{w}{\mu(B)}
 \le \alpha_0w\log\frac{w}{\mu(A_0)}
     +\alpha_1w\log\frac{w}{\mu(A_1)}.
\]
These are the entropy contributions obtained by making the density constant
relative to \(\mu\) on each cell. We will sum the displayed inequality over
matched cells and then pass from cells to points.

If \(\mu\) is a point mass, the result is immediate. Otherwise, work in the
affine hull of its support and enclose the support in a cube \([0,M]^k\).
For each \(P_j\), successively cut each existing cell into \(n\) equal-mass
pieces along coordinates \(1,\ldots,k\), obtaining closed rectangles
\(A_{j,i}\) of mass \(w=n^{-k}\), indexed by the successive choices
\(i=(i_1,\ldots,i_k)\).

These cuts exist because coordinate hyperplanes are null. Each
coordinate \(u\) is nonconstant on the support. Log-concavity applied to
half-spaces, using inner regularity, shows that \(S_u(t)=\mu\{u>t\}\) is
log-concave. An atom at \(a\) with \(S_u(a)>0\) would produce a jump of
\(\log S_u\) inside \(\{S_u>0\}\), contradicting its concavity. Applying the
same argument to \(-u\) excludes an atom at the upper endpoint. Since
\(P_j\ll\mu\), the cut planes are null for \(P_j\) as well.

Pair rectangles with the same index and put
\(B_i=\alpha_0A_{0,i}+\alpha_1A_{1,i}\), which also partition the cube up to
null boundaries. After each common prefix, the cut points for the interpolated
rectangles are weighted averages of two ordered lists, so they remain
ordered and run from \(0\) to \(M\).
Since \(P_j\ll\mu\), each \(\mu(A_{j,i})>0\), and log-concavity then gives
\(\mu(B_i)>0\).

Choose an index \(I\) uniformly and, given \(I=i\), draw \(X_{j,n}\) from
\(P_j(\,\cdot\mid A_{j,i})\) and \(Y_n\) from \(\mu(\,\cdot\mid B_i)\).
Then \(X_{j,n}\sim P_j\), and both
\(Z_n=\alpha_0X_{0,n}+\alpha_1X_{1,n}\) and \(Y_n\) lie in \(B_I\).
The law \(Q_n\) of \(Y_n\) assigns mass \(w\) to each \(B_i\) and has constant density
\(w/\mu(B_i)\) there relative to \(\mu\). Summing the cell inequality and
applying Jensen's inequality for \(t\log t\) to the endpoint densities on
each cell gives
\[
 \KL{Q_n}{\mu}\le\alpha_0\KL{P_0}{\mu}+\alpha_1\KL{P_1}{\mu}.
\]

It remains to show that \(Y_n\) and \(Z_n\) agree in the limit. Given the
preceding cuts, the next coordinate interval is chosen uniformly from \(n\)
intervals whose lengths sum to \(M\), so its mean length is \(M/n\).
Since \(Y_n\) and \(Z_n\) lie in the same rectangle,
\[
 \E\|Y_n-Z_n\|_1\le\frac{kM}{n}.
\]
By compactness, a subsequence of \((X_{0,n},X_{1,n},Y_n)\) converges in law
to some \((X_0,X_1,Y)\). The marginals satisfy \(X_j\sim P_j\), and the
distance estimate gives \(Y=\alpha_0X_0+\alpha_1X_1\) almost surely.
Finally, use the variational formula
\begin{equation}\label{eq:variational}
 \KL P Q=\sup_g\left(\int g\dd P-\log\int e^g\dd Q\right),
\end{equation}
where the supremum is over bounded continuous \(g\); it follows from
\citet[equation~(4.5)]{asz} by optimising over additive constants.
This expresses \(P\mapsto\KL P\mu\) as a supremum of weakly continuous functions,
so the entropy bound for \(Q_n\) passes to \(\law(Y)\), giving the desired result.
\end{proof}

To compare selection gains, fix affine functions \(v_1,\ldots,v_\ell\)
on \(\R^d\) and a common index distribution \(\rho\).
For compactly supported \(\xi\), define the largest expected payoff at budget \(h\) by
\[
 R_\xi(h)=\sup\{\E v_L(W): W\sim\xi,\ L\sim\rho,\ I(L;W)\le h\},
 \qquad h\ge0,
\]
where \(I(L;W)=\KL{\law(L,W)}{\rho\otimes\xi}\) is the mutual information.

\begin{proposition}[Comparison of information budgets]\label{prop:budgets}
Let \(\mu\) be a compactly supported log-concave probability law on \(\R^d\),
and let \(D=\KL\nu\mu<\infty\) for a probability law \(\nu\).
For \(a,c\ge0\) and \(0<\theta<1\), put \(\eta=(1-\theta)a+\theta(c+D)\).
Then
\begin{equation}\label{eq:common}
 (1-\theta)R_\mu(a)+\theta R_\nu(c)
 \le (1-\theta)R_\mu(\eta)+\theta R_\nu(\eta).
\end{equation}
Consequently, for \(h\ge0\) and \(\delta>0\),
\begin{equation}\label{eq:gains}
 R_\nu(h+\delta)-R_\nu(\delta)
 \le \left(1+\frac{D}{\delta}\right)\bigl(R_\mu(h)-R_\mu(0)\bigr).
\end{equation}
\end{proposition}

\begin{proof}
Take couplings at budgets \(a\) under \(\mu\) and \(c\) under \(\nu\), using
a common label \(L\sim\rho\). For any such pair \((L,W)\) with \(W\sim\xi\),
the conditional laws \(P_i=\law(W\mid L=i)\) satisfy the entropy chain rule
\begin{equation}\label{eq:chain}
 \sum_{i:\rho_i>0}\rho_i\KL{P_i}{\mu}
 =\KL{\law(L,W)}{\rho\otimes\mu}=I(L;W)+\KL{\xi}{\mu}.
\end{equation}
Indeed, \(\dd P_i/\dd\mu=(\dd P_i/\dd\xi)(\dd\xi/\dd\mu)\) almost surely
under \(P_i\). For each label with \(\rho_i>0\), both conditional endpoint laws
therefore have finite entropy, so we couple them by Lemma~\ref{lem:interp}.
Writing \(X_0\sim\mu\), \(X_1\sim\nu\) and \(Z=(1-\theta)X_0+\theta X_1\),
we average \eqref{eq:interp} over the labels and apply \eqref{eq:chain} to obtain
\[
 I(L;Z)+\KL{\law(Z)}{\mu}\le (1-\theta)a+\theta(c+D)=\eta.
\]

Resampling preserves the expected payoff at the interpolated point and gives a label
whose mutual information with either endpoint is at most \(\eta\).
Draw \(L'\) from \(\law(L\mid Z)\),
independently of the endpoints given \(Z\),
as in the proof of \citet[Lemma~3.2]{hww}. Since \((L',Z)\) has the same law as \((L,Z)\),
data processing along \(X_j\to Z\to L'\) gives
\(I(L';X_j)\le I(L';Z)=I(L;Z)\le\eta\). By affinity,
\[
 \begin{aligned}
 (1-\theta)\E v_L(X_0)+\theta\E v_L(X_1)
   &=\E v_L(Z)=\E v_{L'}(Z)\\
   &\le (1-\theta)R_\mu(\eta)+\theta R_\nu(\eta).
 \end{aligned}
\]
Taking the suprema over the two original couplings proves \eqref{eq:common}.

To deduce \eqref{eq:gains}, note that \(R_\nu\) is concave because, with the
marginals fixed, the payoff is linear in the coupling and mutual information is convex.
For \(x\ge0\) and \(0<t\le\delta\), take \(a=x\) and \(c=x+t+\delta\).
We choose the common budget \(\eta=x+t\) so that the increments from \(a\)
to \(\eta\) and from \(\eta\) to \(c\) have lengths \(t\) and \(\delta\),
respectively, which determines \(\theta=t/(t+\delta+D)\).
Comparing secant slopes gives
\[
 \frac{R_\nu(c)-R_\nu(c-t)}{t}
 \le \frac{R_\nu(c)-R_\nu(c-\delta)}{\delta}
 \le \left(1+\frac{D}{\delta}\right)\frac{R_\mu(x+t)-R_\mu(x)}{t}.
\]
The first inequality is concavity on \([c-\delta,c]\), and the second is
\eqref{eq:common}. Since \(c-t=x+\delta\), multiplying by \(t\) and summing
over a partition of \([0,h]\) of mesh at most \(\delta\) gives \eqref{eq:gains}.
\end{proof}

\begin{proposition}[Laplace comparison]\label{prop:laplace}
Let \(X\in\R^d\) have log-concave law \(\mu\), and let \(N\) be a seminorm.
Fix \(\lambda>0\) and put \(w=\E N(X-\E X)\) and
\(K_\lambda=\sup_{u\in\bb}\E e^{\lambda u(X)}\).
For \(r>1\),
\begin{equation}\label{eq:laplace}
 \log\E e^{\lambda N(X)/r}
 \le \frac{\lambda w}{r-1}+\frac1r\log K_\lambda.
\end{equation}
Equivalently, for every probability law \(\nu\) with
\(D=\KL\nu\mu<\infty\),
\begin{equation}\label{eq:entropyform}
 \lambda\int N\dd\nu\le\log K_\lambda+
 \bigl(\sqrt D+\sqrt{\lambda w}\bigr)^2.
\end{equation}
\end{proposition}

\begin{proof}
\begin{samepage}
Fix \(r>1\). Applying the Gibbs variational formula to \(g=\lambda N/r\),
by truncation since \(N\ge0\), gives
\[
 \log\E_\mu e^{\lambda N(X)/r}
 =\sup_{\nu:\,\KL\nu\mu<\infty}
 \left\{\frac{\lambda}{r}\E_\nu N-\KL\nu\mu\right\}.
\]
To prove \eqref{eq:laplace}, we need to bound the expression in braces uniformly
over \(\nu\), so fix such a law and write \(D=\KL\nu\mu\).
\end{samepage}

Suppose first that \(\mu\) has compact support. Let \(f=\max_i u_i\) for
\(u_1,\ldots,u_\ell\in\bb\). Under \(\nu\), the least index attaining \(f\) has
law \(\rho\) and mutual information
\(H(\rho)=-\sum_{\rho_i>0}\rho_i\log\rho_i\).
For \(h=H(\rho)\) and \(\delta>0\), \eqref{eq:gains} therefore gives
\[
 \E_\nu f=R_\nu(h+\delta)
 \le R_\nu(\delta)+\left(1+\frac D\delta\right)
 \bigl(R_\mu(h)-R_\mu(0)\bigr).
\]
Under \(\mu\), independence gives \(R_\mu(0)=\sum_i\rho_i u_i(\E X)\),
so \(\E u_L(X)-R_\mu(0)=\E u_L(X-\E X)\le w\) for any selection
\(L\sim\rho\). Thus \(R_\mu(h)-R_\mu(0)\le w\).

Under \(\nu\), a rule at budget \(\delta\) has joint entropy
at most \(D+\delta\) relative to \(\rho\otimes\mu\), by \eqref{eq:chain}.
With \(A=\lambda w\), testing \eqref{eq:variational} with
\((i,x)\mapsto\lambda u_i(x)\) and combining these bounds gives
\[
 \begin{aligned}
 \lambda R_\nu(\delta)
 &\le D+\delta+\log\!\left(\sum_i\rho_i\E_\mu e^{\lambda u_i(X)}\right)
 \le D+\delta+\log K_\lambda,\\
 \lambda\E_\nu f&\le \log K_\lambda+D+A+\delta+AD/\delta.
 \end{aligned}
\]
For \(A>0\), choose \(\delta=A/(r-1)\) so that the coefficient of \(D\) is \(r\),
and for \(A=0\), let \(\delta\downarrow0\).
By Hahn--Banach, we can take nonnegative finite maxima increasing to \(N\),
so monotone convergence gives the same bound for \(\E_\nu N\) and the opening
variational formula proves \eqref{eq:laplace}.

To remove compact support, first note that each unbounded upper tail decays
exponentially, since \(\log S_u\) lies below the extension
of a decreasing secant. Applying this to both signs of the coordinates
gives moments of all orders for \(X\) and \(N(X)\).
Condition on the closed balls \(B_j\) of radius \(j\) about zero, of mass
\(c_j\uparrow1\). For \(c_j>0\), the conditional law \(\mu_j\) is log-concave
by convexity of \(B_j\).
The conditional means converge to \(\E X\), and
\(|N(x-a)-N(x-b)|\le N(a-b)\) gives
\(w_j=\E_{\mu_j}N(X-\E_{\mu_j}X)\to w\).
In \eqref{eq:laplace} for \(\mu_j\), we have \(K_{\lambda,j}\le K_\lambda/c_j\)
and \(c_j\E_{\mu_j}e^{\lambda N(X)/r}\uparrow\E_\mu e^{\lambda N(X)/r}\)
by monotone convergence. Taking limits proves \eqref{eq:laplace}.
Finally, Gibbs duality and
\(\inf_{r>1}\{rD+rA/(r-1)\}=(\sqrt D+\sqrt A)^2\) give \eqref{eq:entropyform}.
\end{proof}

\begin{corollary}\label{cor:tail}
Let \(X\in\R^d\) be log-concave and \(N\) a seminorm. For \(t\ge16m\),
\[
 \Prob(N(X)>t)\le\sup_{u\in\bb}\Prob(u(X)>t/8),
\]
with strict inequality whenever \(\Prob(N(X)>t)>0\).
\end{corollary}
\begin{proof}
Assume \(q=\Prob(N(X)>t)>0\), and put \(h=\log(1/q)\) and \(a=t/8\).
Then \(t\ge16m>0\) and Markov's inequality gives \(h\ge\log16\).
For a contradiction, suppose \(S_u(a)\le q\) for all \(u\in\bb\).
We will show that \(\log K_\lambda\le11h/16\) at \(\lambda=7h/(2t)\).
For \(\nu=\mu(\,\cdot\mid N>t)\), we have \(\KL\nu\mu=h\) and
\(\E_\nu N>t\). Since \(w\le2m\le t/8\), this bound and
\eqref{eq:entropyform} would give the contradiction
\[
 \frac72h<\lambda\E_\nu N
 \le\frac{11h}{16}+\left(\sqrt h+\sqrt{\frac{7h}{16}}\right)^2
 =\left(\frac{17}{8}+\frac{\sqrt7}{2}\right)h<\frac72h.
\]
To establish the bound, note that symmetry of \(\bb\) gives \(S_u(-a)\ge1-q\),
since nonconstant projections are atomless and constant projections lie in \((-a,a)\).
Log-concavity gives \(S_u(x)\le q e^{-\kappa(x-a)}\) for \(x\ge a\), where
\(\kappa=4(h+\log(1-q))/t\ge15h/(4t)\), since
\(-\log(1-q)\le q/(1-q)\le1/15\le h/16\).
Integrating the tail at \(\lambda\) gives
\[
 \log K_\lambda\le\lambda a+\log\left(1+\frac{q\lambda}{\kappa-\lambda}\right)
 \le\frac{7h}{16}+\log(15/8)\le\frac{11h}{16}. \qedhere
\]
\end{proof}

For independent centred coordinates with regular moment growth,
\citet[Corollary~1.3]{lsind} proved a related tail comparison in terms of the
two-sided tails of individual linear functionals.

\section{Proof of Theorem~\ref{thm:main}}

The Laplace comparison gives, for \(1\le p<\infty\),
\begin{equation}\label{eq:momentlaplace}
 \|N(X)\|_p\le w+
 \inf_{\lambda>0:\,K_\lambda<\infty}\frac{p+\log K_\lambda}{\lambda}.
\end{equation}
To obtain \eqref{eq:momentlaplace}, take \(r=1+(\lambda w+\log K_\lambda)/p>1\) for \(K_\lambda<\infty\).
Since \(0\in\bb\), we have \(\log K_\lambda\ge0\), so
\eqref{eq:laplace} gives \(\log\E e^{\lambda N(X)/r}\le p\), and
maximising \(y^p e^{-\lambda y/r}\) over \(y\ge0\) gives
\(\|N(X)\|_p\le pr(e\lambda)^{-1}(\E e^{\lambda N(X)/r})^{1/p}
\le pr/\lambda\). If \(\lambda w+\log K_\lambda=0\), let \(r\downarrow1\).

Applying \eqref{eq:momentlaplace} requires bounding \(\E e^{\lambda u(X)}\) at a
common value of \(\lambda\) for all \(u\in\bb\). Each projection \(u(X)\)
is log-concave and satisfies \(\|u(X)\|_p\le s_p\), so a one-dimensional estimate suffices.

\begin{lemma}\label{lem:scalar}
Let \(T\) be a real log-concave random variable with \(\|T\|_p\le s\),
where \(1\le p<\infty\) and \(s>0\). With
\(\lambda=(p+1)/(2es)\), one has \(\log\E e^{\lambda T}\le 1+\lambda s\).
\end{lemma}

\begin{proof}
By scaling, it suffices to prove the result for \(s=1\).
Put \(q=\Prob(T>-1)\). If \(q=0\), then
\(\E e^{\lambda T}\le e^{-\lambda}\), proving the claim.
Otherwise, \(\E|T|^p\le1\) and \(|T|^p\ge1\) on \(\{T\le-1\}\) give
\[
 q\,\E(|T|^p\mid T>-1)\le1-\Prob(T\le-1)=q.
\]
Let \(Y\) have the law of \(1+T\) given \(T>-1\).
The bound \(\E(|T|^p\mid T>-1)\le1\) and Minkowski's inequality give
\(\|Y\|_p\le2\), while \(\Prob(Y>t)=\Prob(T>t-1)/q\) is log-concave for \(t\ge0\).

To control higher moments, we use the comparison for log-concave survival functions
\citep[Proposition~3.3]{flm}. Let \(G_p=\Gamma(p+1)^{1/p}\) be the \(L^p\)-norm
of an exponential random variable of mean \(1\), and choose an exponential \(Z\)
with mean \(\|Y\|_p/G_p\), so that \(\E Z^p=\E Y^p\).
Since \(t\mapsto\log\Prob(Y>t)+t/\E Z\) is concave and vanishes at zero,
\(\Delta(t)=\Prob(Y>t)-\Prob(Z>t)\) can change sign only from positive to negative.
As \(\int_0^\infty t^{p-1}\Delta(t)\dd t=0\), the increasing weight
\(t^{k-p}\) favours the negative part for \(k>p\), giving
\[
 \E Y^k\le\E Z^k=\Gamma(k+1)(\|Y\|_p/G_p)^k\qquad(k\ge p).
\]

To sum these estimates, put \(a=2\lambda/G_p\).
Stirling's lower bound gives \(eG_p\ge p+3/2\), since the function
\(\tfrac12\log(2\pi p)-p\log(1+3/(2p))\) is increasing for \(p\ge1\)
and positive at \(1\). Thus \(a\le(2p+2)/(2p+3)<1\).
Using \(\|Y\|_k\le2\) for \(1\le k\le p\)
and the comparison above for \(k>p\) gives
\[
 \E e^{\lambda Y}
 \le e^{2\lambda}+\frac{a^{\lfloor p\rfloor+1}}{1-a}
 \le e^{2\lambda}+B_p,\qquad
 B_p=(2p+3)\left(\frac{2p+2}{2p+3}\right)^{\max\{p,2\}}.
\]
The function \(B_p e^{-(p+1)/e}\) is increasing on \([1,2]\)
and decreasing on \([2,\infty)\). Its maximum is
\(36/(7e^{3/e})<e-1\), so \(\E e^{\lambda Y}\le e^{2\lambda+1}\).
Conditioning on \(\{T>-1\}\) and using \(\E e^{\lambda Y}\ge1\) gives
\(\E e^{\lambda T}\le e^{-\lambda}\E e^{\lambda Y}\le e^{\lambda+1}\).
\end{proof}

\begin{proof}[Proof of Theorem~\ref{thm:main}]
For a countable dense subset \(\mathcal U\) of \(\bb\), we have \(N(X)=\sup_{u\in\mathcal U}|u(X)|\),
so \(\|N(X)\|_\infty=s_\infty\). This also gives \(N(X)=0\) almost surely if \(s_p=0\).

For \(p<\infty\) and \(s_p>0\), Lemma~\ref{lem:scalar} gives
\(\log K_\lambda\le1+\lambda s_p\) at \(\lambda=(p+1)/(2es_p)\).
Thus \eqref{eq:momentlaplace} yields
\(\|N(X)\|_p\le w+s_p+(p+1)/\lambda=w+(1+2e)s_p\).
Since \(w\le m+N(\E X)\le m+s_p\), this proves \eqref{eq:main}.
\end{proof}

\begin{remark}
Any universal bound \(\|N(X)\|_p\le C_1 m+C s_p\) with \(C_1,C\ge0\) forces \(C_1\ge1\).
For a standard Gaussian in \(\R^d\), the Euclidean norm gives
\(m\le\sqrt d=\|N(X)\|_2\) and \(s_2=1\), so \((1-C_1)\sqrt d\le C\) in all dimensions.

With the coefficient of \(m\) fixed at one, let \(C_*\) be the least universal
coefficient of \(s_p\). Theorem~\ref{thm:main} gives \(C_*\le C_0\).
To show that \(C_*>1\), take the Euclidean norm \(N\) and density
\(e^{-N(x)}/(8\pi)\) on \(\R^3\).
The radius \(R=N(X)\) has density \(r^2e^{-r}/2\) and mean \(3\).
The direction of \(X\) is independent of \(R\) and has first coordinate uniform on \([-1,1]\),
so \(\|R\|_p=(p+1)^{1/p}s_p\) and \(s_p^p=(p+2)\Gamma(p+1)/2\).
Thus \(s_p\ge G_p\ge p/e\), and \(e^x\ge1+x\) gives
\(\|R\|_p-s_p\ge\log(p+1)/e\to\infty\).
Since \(\E R=3\), some finite \(p\) gives
\(C_*\ge(\|R\|_p-3)/s_p>1\).
\end{remark}

\clearpage
\begingroup
\singlespacing
\section*{Acknowledgements}
GPT-6 Astra was used to assist with the literature search and to check and refine the arguments.

\endgroup

\end{document}